\documentclass[10pt]{amsart}
\pdfoutput=1
\usepackage{amsmath, amsthm, amssymb,slashed}
\usepackage{ifpdf}
\usepackage[pdftex]{graphicx}
\usepackage{tikz}
\usetikzlibrary{matrix,arrows,calc}

\usepackage[pdftex,plainpages=false,hypertexnames=false,pdfpagelabels]{hyperref}
 \theoremstyle{plain}
 \newtheorem{thm}{Theorem}[section]

 \newtheorem{lem}[thm]{Lemma}
  
 \newtheorem{prop}[thm]{Proposition}
 \theoremstyle{definition}
 \newtheorem{defn}[thm]{Definition}

 \theoremstyle{remark}
 \newtheorem{rmk}[thm]{Remark}

\def\beq{\begin{eqnarray}}
\def\eeq{\end{eqnarray}}

\DeclareSymbolFont{bbold}{U}{bbold}{m}{n}
\DeclareSymbolFontAlphabet{\mathbbold}{bbold}

 \newcommand{\bp}{\begin{proof}[Proof]}
 \newcommand{\ep}{\end{proof}}
 
\DeclareMathOperator{\SM}{\underline{\sf SMfld}}

\DeclareMathOperator{\Euc}{\underline{\sf Euc}}

\DeclareMathOperator{\pt}{\rm pt}
\DeclareMathOperator{\vac}{{{{\rm vac}}}}

\DeclareMathOperator{\PA}{{{{\rm PA}}}}
\DeclareMathOperator{\ev}{{\rm ev}}

\newcommand{\sq}{\mathord{/\!/}}

\def\R{{\mathbb{R}}}

\def\id{{{\rm id}}}
\def\proj{{\rm proj}}

\def\C{{\mathbb{C}}}

\def\T{{\mathbb{T}}}
\def\Z{{\mathbb{Z}}}

\def\pt{{\rm pt}}

\def\SM{ {\underline{\sf SMfld}}}

\def\dashdownarrow{\ensuremath{\rotatebox[origin=c]{-90}{$\dashrightarrow$}}}

\begin{document}

\title[Perturbative supersymmetric quantum mechanics and ${\rm L}$-theory]{Perturbative $N=2$ supersymmetric quantum mechanics and ${\rm L}$-theory with complex coefficients}

\author{Daniel Berwick-Evans}

\address{Department of Mathematics, Stanford University}

\email{danbe@stanford.edu}

\date{\today}

\maketitle 

\begin{abstract}
We construct L-theory with complex coefficients from the geometry of $1|2$-dimensional perturbative mechanics. Methods of perturbative quantization lead to wrong-way maps that we identify with those coming from the MSO-orientation of L-theory tensored with the complex numbers.
\end{abstract}

\setcounter{tocdepth}{1}
\tableofcontents


\section{Introduction}

In this paper we begin an investigation of $d|2$-dimensional super Euclidean field theories (EFTs) in the style of Stolz-Teichner, starting with~$d=1$. Since $1|2$-dimensional quantum mechanics has twice the minimal number of fermions required for supersymmetry, it often goes by the moniker \emph{$N=2$ supersymmetric quantum mechanics}. Concrete examples are provided by the Hodge-de~Rham complex of an oriented Riemannian manifold, as explained by Witten in his study of Morse theory~\cite{susymorse} and reviewed briefly in Section~\ref{sec:background} below.

The exact geometrical and topological information captured by $1|2$-dimensional theories remains unclear to us, and we view the present paper as extracting the most easily computable quantities. In Theorem~\ref{thm1} we show that the fiberwise partition function for a family of field theories over~$X$ defines a class in~${\rm L}(X)\otimes \C$, where~${\rm L}$ is the cohomology theory defined by the symmetric ${\rm L}$-theory spectrum. We note that---being a cohomology theory over~$\C$---${\rm L}\otimes \C$ is nothing more than 4-periodic de~Rham cohomology. Our reason for viewing it in this more complicated light is the natural appearance of Hirzebruch's ${\rm L}$-genus: a perturbative quantization procedure constructs a pushforward to the point that agrees with one coming from the MSO-orientation of L-theory tensored with~$\C$ (see Theorem~\ref{thm2}). 
In particular, $\pi_!\colon {\rm L}^0(X)\otimes \C\to {\rm L}^{-n}(\pt)\otimes \C$ sends~1 to the ${\rm L}$-genus (or signature) of~$X$. 

Our long-term goal is to leverage an understanding of the $1|2$-dimensional case to gain traction on the more complicated $2|2$- and $3|2$-dimensional theories. This is in analogy with Stolz and Teichner's approach to a geometric model for elliptic cohomology: they are motivated in large part by the relation between $1|1$-dimensional EFTs and Dirac operators on Riemannian spin manifolds~\cite{mingeo, ST11}. In the footsteps of G.~Segal~\cite{Segal_Elliptic} they argue by analogy that $2|1$-EFTs ought to capture structures related to Dirac operators on \emph{loop spaces}. 
Nuances of super Euclidean geometry rule out EFTs of dimension $d|1$ for~$d>2$ (e.g., see~\cite{5lectures}), so the next class of Euclidean field theories to consider are of dimension $d|2$, which exist for~$d=0,1,2,3$. Unbridled optimism might lead one to hope for a relationship between $2|2$-dimensional field theories, TMF with level structure, and the Ochanine genus. However, as yet such connections are purely by analogy to the lower dimensional case described below, coupled to Witten's description of the Ochanine genus as the signature of loop space~\cite{Witten_Dirac}.

\subsection{Statement of results}

Drawing on methodology from \cite{DBE_WG}, we extract homotopical data from $1|2$-Euclidean field theories associated to the moduli of $1|2$-dimensional super circles. We focus on a connected component of this moduli space characterized by the holonomy of a~$0|2$-dimensional odd vector bundle over an ordinary circle which takes the form ${\rm diag}(+1,-1)$; we call these \emph{periodic-antiperiodic}~(PA) circles. For a smooth manifold~$X$ we define a stack~$\Phi^{1|2}_{0,\PA}(X)$ that consists of periodic-antiperiodic circles with a map to~$X$ of zero classical action. There is a line bundle over~$\Phi_{0,\PA}^{1|2}(X)$ denoted $\mathcal{L}$. We call a section of~$\mathcal{L}^{\otimes k}$ \emph{supersymmetric} if it extends to a section over a stack~$\Phi^{1|2}_{\epsilon,\PA}(X)$, that roughly consists of PA super circles mapping to~$X$ with \emph{nilpotent} action. Tensor products of line bundles yields a graded super commutative algebra of supersymmetric sections. Our first result is the following. 

\begin{thm} \label{thm1}
Concordance classes of supersymmetric sections of $\mathcal{L}^\bullet$ over the periodic-antiperiodic vacua of~$X$ are naturally isomorphic to 4-periodic de~Rham cohomology,
$$
\bigoplus_{i\in \Z} {\rm H}_{\rm dR}^{\bullet+4i}(X) \cong \Gamma_{\rm susy}(\Phi_{0,\PA}^{1|2}(X);\mathcal{L}^\bullet)/{\rm concordance}
$$
where the cup product on the left coincides with the tensor product of sections on the right. 
\end{thm}

The symmetric ${\rm L}$-theory spectrum defines a 4-periodic cohomology theory that after tensoring with~$\C$ is 4-periodic de~Rham cohomology, ${\rm L}^k(X)\otimes \C\cong \bigoplus_{i\in \Z} {\rm H}_{\rm dR}^{k+4i}(X)$. Hence, we can rephrase the above result in terms of an isomorphism with ${\rm L}$-theory with complex coefficients. The connection between $1|2$-dimensional field theories and ${\rm L}$-theory does not end here: our second main result is a geometric construction of wrong-way maps
\begin{equation}
\begin{array}{c}
\begin{tikzpicture}[node distance=3.5cm,auto]
  \node (A) {$\Gamma_{\rm susy}(\Phi_{0,\PA}^{1|2}(X);\mathcal{L}^k)$};
  \node (B) [node distance= 4.5cm, right of=A] {$\Gamma_{\rm susy}(\Phi_{0,\PA}^{1|2}(\pt);\mathcal{L}^{k-n})$};
  \node (C) [node distance = 1.5cm, below of=A] {${\rm L}^k(X)\otimes \C$};
  \node (D) [node distance = 1.5cm, below of=B] {${\rm L}^{k-n}(\pt)\otimes \C$};
  \node (E) [node distance =.75cm, below of=A] {$\null$};
    \node (F) [node distance =2.25cm, right of=E] {$\circlearrowright$};
  \draw[->] (A) to node {$\pi_!$} (B);
  \draw[->>] (A) to (C);
  \draw[->] (C) to node {$\pi_!$} (D);
  \draw[->>] (B) to (D);
\end{tikzpicture}\end{array}
\nonumber \quad {\rm dim}(X)=n
\end{equation}
compatible with those coming from the ${\rm MSO}$-orientation of ${\rm L}$-theory tensored with~$\C$.

This construction starts with the linearization of the fields of the classical $1|2$-sigma model over~$\Phi_{0,\PA}^{1|2}(X)$, which defines an infinite dimensional vector bundle over~$\Phi^{1|2}_{0,\PA}(X)$. The linearization of the classical action defines a family of operators on the fibers denoted by~$\Delta_{X,\PA}^{1|2}$, called the \emph{kinetic operators}. We define the $\zeta$-super determinant of the kinetic operators as
$$
{\rm sdet}_\zeta(\Delta_{X,\PA}^{1|2}):=\frac{{\rm pf}_\zeta(\Delta_{X,\PA}^{1|2}|_{{\rm fer}})}{{\rm det}_\zeta(\Delta_{X,\PA}^{1|2}|_{{\rm bos}})^{1/2}}
$$
where the numerator is the $\zeta$-regularized Pfaffian on odd (or fermionic) sections and the denominator is the square root of the $\zeta$-determinant on even (or bosonic) sections.  Our convention is to take the positive square root of the determinant. From our point of view, perturbative quantization is the construction of a volume form on the fibers~$\Phi^{1|2}_{0,\PA}(X)\to \Phi^{1|2}_{0,\PA}(\pt)$ which in turn defines a wrong-way map. Explicitly, such a volume form is the product of the super determinant of the kinetic operator with the canonical volume form that comes from integration of differential forms.

\begin{thm} \label{thm2} The $\zeta$-super determinant of $\Delta_{X,\PA}^{1|2}$ represents the ${\rm L}$-class of $X$ as a supersymmetric function on $\Phi_{0,\PA}^{1|2}(X)$. This function modifies the canonical volume form on~$\Phi_{0,\PA}^{1|2}(X)$ so that integration of supersymmetric sections of $\mathcal{L}^\bullet$ coincides with the wrong-way map associated with the ${\rm MSO}$-orientation for ${\rm L}$-theory with complex coefficients. In particular, the total volume of~$\Phi_{0,\PA}^{1|2}(X)$ is the ${\rm L}$-genus of~$X$. 
\end{thm}

The $\zeta$-determinant computation underlying Theorem~\ref{thm2} is of a fairly standard variety in the physics literature, so in this regard our contribution is to iron out some mathematical details, explain the geometry leading to the relevant operators, and contextualize the result via Theorem~\ref{thm1}.


\begin{rmk}
The analogous constructions for periodic-periodic super circles (i.e., those defined by an $\R^{0|2}$-bundle over $S^1$ with trivial holonomy) give quite different results. Owing to a failure of excision, supersymmetric sections for tensor powers of any line bundle over periodic-periodic vacua cannot be given the structure of a nontrivial cohomology theory. However, the only proof we know of this fact is by brute force, and therefore rather unenlightening (see~\cite{DBEPhd}, Chapter~4). A computation similar to the one constructing the L-class shows that the $\zeta$-super determinant of the relevant kinetic operator over periodic-periodic circles~1. However, in this case the classical action defines an interesting volume form on action zero fields that we analyzed in depth in~\cite{DBE_CGB}. With respect to this volume form, the total volume of periodic-periodic vacua is the Euler characteristic of~$X$. In physics parlance, the contribution to the path integral of the nonzero modes is~1, and the contribution from the zero modes is the Euler characteristic. This result can be anticipated from the physical interpretation of the index theorem~\cite{Alvarez}:~$N=2$ supersymmetric quantum mechanics of an oriented Riemannian manifold~$X$ encodes two super traces of the operator $d+d^*$ acting on $\Omega^\bullet(X)$, one corresponding to the signature, and the other corresponding to the Euler characteristic, e.g., see~\cite{susymorse}. 
\end{rmk}

\subsection{Notation and conventions} This paper follows the conventions for supermanifolds and super stacks set forth in~\cite{DBE_WG}, and we refer the reader to the subsection therein with the same title as the present one for details. In summary, our supermanifolds have complex algebras of functions (called $cs$-manifolds in~\cite{DM}), we work with the functor of points throughout, and we frequently identify a stack with a groupoid presenting it as articulated in~\cite{BlohmannStacks}. One important convention for group actions on mapping spaces (or stacks) is that the action of $\phi\in {\rm Aut}(\Sigma)$ on $f\in {\sf Map}(\Sigma,X)$ is by~$f\mapsto f\circ \phi^{-1}$. 

\subsection{Acknowledgements} Part of this work was completed while I was a Ph.D. student, and I warmly thank my advisor Peter Teichner for his support and guidance. I also thank Owen Gwilliam for several helpful conversations in the early stages of this work. 

%
%
%
%

\section{Background: $1|2$-dimensional geometry and physics}\label{sec:background}

\subsection{$1|2$-Euclidean manifolds}\label{sec:modelgeo}
Let~$\R^{1|2}$ denote the super group with multiplication
$$
(t,\theta_1,\theta_2)\cdot (t',\theta_1',\theta_2')=(t+t'+i\theta_1\theta_1'+i\theta_2\theta_2',\theta_1+\theta_1',\theta_2+\theta_2'), \quad (t,\theta_1,\theta_2),(t',\theta_1',\theta_2')\in \R^{1|2}(S).
$$
This defines the $1|2$-dimensional \emph{super translation group}. The Lie algebra of right-invariant vector fields on this super Lie group has generators $D_1:=\partial_{\theta_1}-i\theta_1\partial_t$ and $D_2:=\partial_{\theta_2}-i\theta_2\partial_t$ satisfying the relations 
$$
D_1^2=\frac{1}{2}[D_1,D_1]=D_2^2=\frac{1}{2}[D_2,D_2]=-i\partial_t, \quad [D_1,D_2]=0.
$$ 
The left-invariant vector fields differ by a sign, e.g.,~$\partial_{\theta_1}+i\theta_1\partial_t$. 

We will be interested in field theories that are invariant under the automorphism that when restricted to~$\R$ is the orientation reversing map~$t\mapsto -t$. Physically, this is the automorphism of \emph{time-reversal}. There are essentially two independent lifts of this automorphism to~$\R^{1|2}$, 
$$
(t,\theta_1,\theta_2)\stackrel{{\sf r}_+}{\mapsto} (t,i\theta_1,i\theta_2),\quad (t,\theta_1,\theta_2)\stackrel{{\sf r}_-}{\mapsto} (t,-i\theta_1,i\theta_2).
$$
We observe that ${\sf r}_+$ and ${\sf r}_-$ both have order~4, commute with each other, and satisfy ${\sf r}_+^2={\sf r}_-^2$. They generate a discrete abelian group we denote by~$T$ which acts on~$\R^{1|2}$ by automorphisms. 

\begin{rmk} Typically time-preserving and time-reversing automorphisms are of different flavors: the former define linear operators on the Hilbert space of states, whereas the latter define \emph{anti}-linear operators, e.g., see Freed and Moore~\cite{Freed_Moore} for an extensive discussion. In this paper we are only considering the value of theories on circles, which requires \emph{invariance} under symmetries rather than equivariance. This allows us to effectively ignore the issue. \end{rmk}

We define the \emph{$1|2$-Euclidean group} as the semidirect product,~$\R^{1|2}\rtimes T$. The left action of $\R^{1|2}\rtimes T$ on $\R^{1|2}$ defines a \emph{$1|2$-Euclidean model space}; see~\cite{mingeo}, Section~6.3 for detail. Such a model space in turn defines a fibered category of super Euclidean manifolds: for any supermanifold $S$ there is a category of \emph{$1|2$-Euclidean manifolds over $S$}, i.e., $S$-families of supermanifolds that admit atlases where each chart is isomorphic to an open submanifold of $\R^{1|2}$ and transition functions are given by restrictions of the action of~$\R^{1|2}\rtimes T$ on $\R^{1|2}$. Between such families one can define \emph{isometries over $S$}, which locally are given by the action of~$\R^{1|2}\rtimes T$. Let $T_+<T$ denote the subgroup that acts by orientation preserving automorphisms of $\R\subset \R^{1|2}$. The subgroup $\R^{1|2}\rtimes T_+<\R^{1|2}\rtimes T$ is the subgroup of \emph{time preserving isometries.} 

Our main examples of $1|2$-Euclidean manifolds will be \emph{$1|2$-super circles}. Let $\R^{1|2}_{>0}$ be the super semigroup given by restricting the structure sheaf and group structure of $\R^{1|2}$ to $\R_{>0}\subset \R$. Given a generator $R\in \R^{1|2}_{>0}(S)$ and $A\in T_+$ we define an $S$-family of super circles as the quotient, $S\times_{(R,A)}\R^{1|2}:=S\times \R^{1|2}/((R,A)\cdot \Z)$, by the right $\Z$-action generated by the composition of the~$R$ and~$A$ actions on $S\times \R^{1|2}$. We view~$R$ as a super radius and~$A$ as a holonomy. When~$A=\id$ (respectively, $A={\sf r}_+{\sf r}_-$) we call the associated circle \emph{periodic-periodic} or \emph{PP} (respectively, \emph{periodic-antiperiodic} or \emph{PA}).

\subsection{A brief tour of $N=2$ supersymmetric mechanics} \label{sec:mech} Fields of the $1|2$-dimensional sigma model consist of maps $\phi\colon S\times I^{1|2}\to X$ for $S\times I^{1|2}$ a proper $S$-family of $1|2$-dimensional Euclidean intervals and~$X$ a fixed Riemannian manifold. These fields have automorphisms coming from super Euclidean isometries~$S\times I^{1|2}\to S\times{I'}^{1|2}$ over~$S$. There is a function on fields called the \emph{classical action}, 
$$
\mathcal{S}(\phi):=\int_{I^{1|2}\times S/S}\langle D_1\phi,D_2\phi\rangle,
$$ 
where the integral is fiberwise over $S$. Being built from right-invariant vector fields, $\mathcal{S}$ is invariant under the left action of super Euclidean isometries. To unpack this action in terms of more familiar geometric quantities, the map~$\phi$ determines a map $x\colon S\times I\to X$ from a family of ordinary 1-manifolds to $X$, and sections $\psi_1,\psi_2\in x^*\pi TX$ of the odd tangent bundle pulled back to this family. Then (after discarding of terms without derivatives) we obtain the action 
\beq
\mathcal{S}(x,\psi_1,\psi_2)=\int_{I\times S/S} \Big(\langle \dot{x},\dot{x}\rangle+\langle\psi_1,\nabla_{\dot{x}}\psi_1\rangle+\langle\psi_2,\nabla_{\dot{x}}\psi_2\rangle+R(\psi_1,\psi_2,\psi_1,\psi_2)\Big)dt\label{eq:action}
\eeq
where $\dot{x}$ denotes the derivative of the path determined by~$x$, $\nabla_{\dot{x}}$ is the covariant derivative along this path for the Levi-Civita connection, and $R$ denotes the curvature of the Levi-Civita connection. 

One can quantize the classical theory using techniques of geometric quantization. The result is a representation of $\R^{1|2}_{>0}$ on the vector space $\Omega^\bullet(X)$ of differential forms on~$X$, which can be viewed as a super time evolution. One mathematical perspective on this time evolution is to package the quantum theory as a functor from a category $1|2$-Euclidean bordisms to the category of (topological) vector spaces; then the semigroup of super intervals acts on this vector space, which is exactly the representation of~$\R^{1|2}_{>0}$, where a super interval has a length in~$\R^{1|2}_{>0}$. Explicitly, the action is determined by the Lie algebra action where~$D_1$ acts by~$d+d^*$ and~$D_2$ acts by~$i(d-d^*)$ on forms, and it is important that both of these operators square to the Hodge Laplacian. We also get an action of~$T_+$ on differential forms that is determined by ${\sf r}_+^2={\sf r}_-^2$ acting as the mod 2 grading involution on $\Omega^\bullet(X)$, and the action of ${\sf r}_+{\sf r}_-$ is by the Hodge star,~$*$, normalized such that~$*^2=+1$. 

As usual in the bordism-type approach, the value of a field theory on a circle is the (super) trace of an operator. In the case above, the value on a periodic-periodic circle is the Euler characteristic of~$X$, and the value on periodic-antiperiodic circle is the signature of~$X$. 

\begin{rmk} In this example, one can extend the action of $T_+$ to one by $O(2)$, where the $SO(2)<O(2)$ action encodes the $\Z$-grading on forms. This leads to a more restrictive notion of field theory, but might end up being the mathematically desirable one in the end. We have chosen to stay with the more general class here, since the restricted versions will still map to ${\rm L}\otimes \C$ by our result. \end{rmk}

\section{Classical vacua}

\subsection{Fields and vacua}

\begin{defn} Define the \emph{stack of fields} for the $1|2$-dimensional sigma model with target~$X$, denoted $\Phi^{1|2}(X)$, as the stack associated to the prestack with objects over $S$ consisting of triples $(R,A,\phi)$ where $R\in \R^{1|2}_{>0}(S)$ and $A\in T_+$ determine a family of super circles and $\phi\colon S\times_{R_A}\R^{1|2}\to X$ is a map. Morphisms over~$S$ consist of commuting triangles
\beq
\begin{tikzpicture}[baseline=(basepoint)];
\node (A) at (0,0) {$S\times_{R_A} \R^{1|2}$};
\node (B) at (3,0) {$S\times_{R'_{A'}} \R^{1|2}$};
\node (C) at (1.5,-1.5) {$X$};
\node (D) at (1.5,-.6) {$\circlearrowright$};
\draw[->] (A) to node [above=1pt] {$\cong$} (B);
\draw[->] (A) to node [left=1pt]{$\phi$} (C);
\draw[->] (B) to node [right=1pt]{$\phi'$} (C);
\path (0,-.75) coordinate (basepoint);
\end{tikzpicture}\label{diag:tricom}
\eeq
where the horizontal arrow is an isomorphism of $S$-families of super Euclidean $1|2$-manifolds. \label{def:fields}
\end{defn}

We will define a substack of $\Phi^{1|2}(X)$ for families of periodic-antiperiodic circles whose map to~$X$ has zero classical action. First, for ${\sf r}_+{\sf r}_-\in T_+$ and $R\in \R_{>0}(S)\subset \R^{1|2}(S)$, hereafter let $S\times_R\R^{1|2}$ denote the family $S\times_{(R,{\sf r}_+{\sf r}_-)}\R^{1|2}$. Next, let $\proj\colon S\times_R\R^{1|2} \to S\times \R^{0|1}$ denote the ${\sf r}_+{\sf r}_-$-invariant projection $\R^{1|2}\to \R^{0|1}$ given by the formula $(t,\theta_1,\theta_2)\mapsto \theta_1$ in the notation of subsection~\ref{sec:modelgeo}. 

\begin{defn} Let $\Phi^{1|2}_{0,\PA}(X)$ be the full substack of $\Phi^{1|2}(X)$ generated by objects over~$S$ defined by $(R,\phi)$ where $R\in \R_{>0}(S)\subset \R^{1|2}(S)$ and $\phi\colon S\times_{R}\R^{1|2}\to X$ is a map that factors through $\proj$. 
\end{defn}

\begin{rmk} A quick check against Equation~\ref{eq:action} verifies that the maps $\Phi^{1|2}_{0,\PA}(X)$ do indeed have zero classical action. In the component description, for a field to have zero action the map $x\colon S\times S^1\to X$ must be constant along the fibers, i.e., must factor through the projection $S\times S^1\to S$. Then, the odd sections~$\psi_1,\psi_2$ must be parallel along the circle, which (given the constancy of $x$) means $\psi_1$ and $\psi_2$ are constant sections. Finally, given the antiperiodic boundary conditions on $\psi_2$, we require~$\psi_2=0$. \end{rmk}


\subsection{Periodic-antiperiodic circles with nilpotent action}

To define supersymmetric sections in the statement of Theorem~\ref{thm1}, we need to define a stack $\Phi^{1|2}_{\epsilon,{\rm PA}}(X)\supset \Phi^{1|2}_{0,{\rm PA}}(X)$ that consists of fields with \emph{nilpotent} classical action. 

Let $R\in \R^{1|1}_{>0}(S)\hookrightarrow \R_{>0}^{1|2}(S)$ for the inclusion $(t,\rho_1)\mapsto (t,\rho_1,0)$. Again, to simplify the notation we let $S\times_R\R^{1|2}:=S\times_{(R,{\sf r}_+{\sf r}_-)}\R^{1|2}$. We shall define morphisms ${\rm proj}_R\colon S\times_{R}\R^{1|2}\to S\times \R^{0|1}$. We observe that the projection~$\R^{1|2}\to \R^{0|1}$ of the previous section is not invariant under the $\Z$-action associated to the lattice $R=(r,\rho_1,0)$ for $\rho_1\ne 0$, so will not define such a map $\proj_R$. Instead, for $(r,\rho_1,0)\in (\R_{>0} \times \R^{0|2})(S)$, we claim the map 
$$
\proj_R\colon \R^{1|2}\times S \to \R^{0|1}\times S,\quad \proj_R(t,\theta_1,\theta_2,s)=(\theta_1-\rho_1\frac{t}{r},s)
$$
is invariant under the left action of $\Z$ on the family $\R^{1|2}\times S$. Denoting the action of the generator by $\mu_R$, we compute
$$
(\proj_R\circ \mu_R)(t,\theta_1,\theta_2,s)=\proj_R(t+r+i\theta_1\rho_1,\theta_1+\rho_1,-\theta_2,s)=\theta_1+\rho_1-\rho_1\frac{t+r+i\rho_1\theta_1}{r}=\theta_1-\rho_1\frac{t}{r}.
$$
This shows the map $\proj_R$ is invariant, as claimed, so defines a map which we also denote by~$\proj_R$
$$
\proj_R\colon S\times_{R}\R^{1|2}\to S\times \R^{0|1}. 
$$
where the source circles $ S\times_{R}\R^{1|2}$ have $R\in \R^{1|1}_{>0}(S)\subset \R^{1|2}_{>0}(S)$.

\begin{defn} Define the stack $\Phi^{1|2}_{\epsilon,{\rm PA}}(X)$ as the full substack of $\Phi^{1|2}(X)$ consisting of pairs~$(R,\phi)$ where $R\in \R^{1|1}(S)\hookrightarrow \R^{1|2}(S)$ as above and $\phi\colon S\times_{R}\R^{1|2}\to X$ is a map that factors through~$\proj_R$. \end{defn}

%

\subsection{Groupoid presentations}\label{sec:11groupoidpresent}

The factorization condition on the map~$\phi$ in the definition of~$\Phi_{\epsilon,\PA}^{1|2}(X)$ allows us to identify it with a morphism $\phi_0\colon S\times \R^{0|1}\to X$ for $\phi=\phi_0\circ \proj_R$. We can therefore find a groupoid that presents the stack $\Phi_{\epsilon,\PA}^{1|2}(X)$ whose objects are $(R,\phi_0)\in (\R^{1|1}_{>0}\times \SM(\R^{0|1},X))(S)$. The morphisms $(R,\phi_0)\to (R',\phi_0')$, are determined by $S$-families of super Euclidean isometries; lifting to the universal cover of the family of super circles, these isometries are determined (non-uniquely) by an $S$-point $(u,\nu_1,\nu_2,A)\in \Euc(\R^{1|2})(S)$. However, not all super Euclidean isometries of~$\R^{1|2}$ descend to ones on super circles, so the first task is to compute those that do. 

We define the group $\R^{1|1}\rtimes T$ by the multiplication
$$
(u,\nu)\cdot(u',\nu')=(u+u'+\nu\nu',\nu+\nu'), \quad (u,\nu),(u',\nu')\in \R^{1|1}(S),
$$
and where generators of~$T$ act as ${\sf r}_+\cdot (u,\nu)=(-u,i\nu)$ and ${\sf r}_-\cdot (u,\nu)=(-u,-i\nu)$. There is an evident inclusion,
\beq
\R^{1|1}\rtimes T\hookrightarrow \R^{1|2}\rtimes T,\quad (u,\nu,{\sf r})\mapsto (u,\nu,0,{\sf r}).\label{eq:incl}
\eeq

\begin{lem}\label{lem:liftaction}
An $S$-point of $\Euc(\R^{1|2})$ descends to map between family of super circles $S\times_{R}\R^{1|2}$ if and only if it is in the image of the inclusion (\ref{eq:incl}).
\end{lem}

\bp First we verify the assertion for super translations. Let $(t,\theta_1,\theta_2)\in R^{1|2}(S)$ denote a lift of coordinates on a family of super circles, let $(u,\nu_1,\nu_2)\in \R^{1|2}(S)$ be a family of super translations, and let~$\mu_R^{\rm PA}$ denote the action of the generator $(R,{\sf r}_+{\sf r}_-)$ of the family of lattices defining the family of super circles, where $R=(r,\rho_1)$. To obtain a well-defined action by super translations, we require 
$$
\mu_R^{\rm PA}\left((u,\nu_1,\nu_2)\cdot (t,\theta_1,\theta_2)\right)= (u,\nu_1,\nu_2)\cdot \mu_R^{\rm PA}(t,\theta_1,\theta_2).
$$
Computing the left side we get
$$
\mu_R^{\rm PA} (u+t+i\nu_1\theta_1+i\nu_2\theta_2,\nu_1+\theta_1,\nu_2+\theta_2)=(u+t+i\nu_1\theta_1+i\nu_2\theta_2+r+i(\nu_1+\theta_1)\rho_1,\nu_1+\theta_1+\rho_1,-(\nu_2+\theta_2)), 
$$
whereas computing the right we get
$$
(u,\nu_1,\nu_2)\cdot (t+r+i\theta_1\rho_1,\theta_1+\rho_1,-\theta_2)=(u+t+r+i\theta_1\rho_1+i\nu_1(\theta_1+\rho_1)-i\nu_2\theta_2,\nu_1+\theta_1+\rho_1,\nu_2-\theta_2)
$$
from which we conclude that the translation action of $(u,\nu_1,\nu_2)\in\R^{1|2}(S)$ descends to an isometry of an $S$-family of super circles if and only if $\nu_2=0$. This action by super translations changes the base-points on supercircles so we obtain a map 
\beq
S\times_{R}\R^{1|2}\to S\times_{R'}\R^{1|2}\label{eq:PAiso}
\eeq
where $R=(r,\rho_1,0)$ and $R'=(r+ 2\nu_1\rho_1,\rho_1,0)\in \R^{1|2}(S)$, which is computed by the conjugation action of a super translation on $R\in \R^{1|2}_{>0}(S)$.

The action of ${\sf r}_\pm$ on $R\in \R^{1|1}(S)$ is $R=(r,\rho_1)\mapsto(-r,\pm i \rho_1)$. For $R\in \R^{1|1}_{>0}(S)$ a generator of a lattice action, the image of the lattice under this map is a lattice with positive generator $R'=(r,\mp i\rho_1)\in \R^{1|1}_{>0}(S)$. The map
$$
(r,\rho_1,t,\theta_1,\theta_2)\mapsto (r,\mp \rho_1,-r,\pm i\theta_1,i\theta_2)
$$
defines an isometry as in (\ref{eq:PAiso}), and so the lemma is proved. 
\ep


The above proof determines the target family of super circles for a fixed source, i.e., we have formulas for $R'$ given $R$ and an isometry. It remains to understand how these isometries affect the maps~$\phi_0\colon S\times \R^{0|1}\to X$. This amounts to finding the arrow that makes the diagram commute
\beq
\begin{array}{ccc} 
\R^{1|2}\times S& \stackrel{{\rm id}\times \proj_R}{\to} & \R^{0|1}\times S\\
 (u,\nu,{\sf r}_\pm)\cdot \downarrow && \dashdownarrow\\
\R^{1|2}\times S& \stackrel{\proj_{R'}}{\to} & \R^{0|1}\times S,
\end{array}\label{diag:dotted}
\eeq
where the left vertical arrow is the action of $\R^{1|1}\rtimes T$ on $\R^{1|2}$.

\begin{lem}
The unique dotted arrow in diagram~(\ref{diag:dotted}) is given by
$$
\theta\mapsto \pm i\left(\nu_1+\theta-\rho_1\frac{u+i\nu_1\theta_1}{r}\right)
$$
for $\theta$ a coordinate on $\R^{0|1}$. 
\end{lem}

\bp
The action by $T$ on $\R^{0|1}$ is determined by the action of generators ${\sf r}_\pm$, and it is easy to verify that the action through~$\pm i$ on $\theta_1$ leads to a commutative square for~$R$ and~$R'$ as in the proof of Lemma~\ref{lem:liftaction}. For super translations, we compute
$$
\proj_{R'}\left((u,\nu_1)\cdot (t,\theta_1,\theta_2)\right)=\proj_{R'}(u+t+i\nu_1\theta_1,\nu_1+\theta_1,\theta_2)=(\nu_1+\theta_1-\rho_1\frac{u+t+i\nu_1\theta_1}{r})
$$
and comparing with $\proj_R(t,\theta_1,\theta_2)$, we deduce the lemma. \ep

Putting this all together, for an $S$-point of $\SM(\R^{0|1},X)\cong \pi TX$ of the form
$$
\phi_0=x+\theta\psi,\quad x\colon C^\infty(X)\to C^\infty(S), \psi\in {\rm Der}_x(C^\infty(X),C^\infty(S))
$$
where $x$ is an algebra homomorphism and $\psi$ is an odd derivation with respect to~$x$, we have an action on $(r,\rho_1,x,\psi)\in (\R^{1|1}\times \pi TX)(S)$ by
\beq
\begin{array}{c}(u,\nu_1)\cdot (r,\rho_1,x,\psi)\mapsto \left(r+ 2i\nu\rho_1,\rho_1,x+ (\nu_1-\frac{\rho_1 u}{r})\psi,(1+\frac{i\rho_1\nu_1}{r})\psi\right).\end{array}\label{eq:PAaction}
\eeq

\begin{defn}
A \emph{surjective map of stacks} is a morphism of stacks that on $S$-points induces an essentially surjective, full morphism of groupoids. 
\end{defn}

We have obtained the following groupoid description. 

\begin{prop}
There is a surjective map of stacks
$$
(\R^{1|1}_{>0}\times \pi TX)\sq (\R^{1|1}\rtimes T)\to \Phi_{\epsilon,\PA}^{1|2}(X),
$$
where the surjection comes from lifting an isometry of super circles to the universal cover. 
\end{prop}

Restricting to $R\in \R_{>0}(S)\subset \R^{1|1}_{>0}(S)$ we obtain the following. 

\begin{prop} \label{prop:grpdpresent}
There is an equivalence of stacks
$$
(\R_{>0}\times \pi TX)\sq (\T^{1|1}_R\rtimes T)\cong \Phi^{1|2}_{0,\PA}(X).
$$
\end{prop}

\begin{rmk}
The above is not quite a quotient stack since the groups~$\T^{1|1}_R\cong \R^{1|1}/R\cdot \Z$ depend on~$R$, and we use the subscript $R$  as a reminder for this subtlety. However, all groups $\T^{1|1}_R$ are isomorphic and in the computations below we will be primarily concerned with invariance under the action of~$\T^{1|1}_R$, so this subtlety gets lost in the wash.  \label{rmk:quot}
\end{rmk}

\section{${\rm L}$-theory with complex coefficients}
\subsection{Line bundles over classical vacua}

Given the groupoid descriptions of the previous section, the homomorphisms
$$
T\cong \Z/4\times \Z/2\stackrel{p_1}{\to} \Z/4\subset \C^\times\cong {\sf Aut}(\C^{0|1})
$$
define odd line bundles denoted $\mathcal{L}$ over $\Phi_{\epsilon,\PA}^{1|2}(X)$ that are natural in~$X$, where the first isomorphism sends ${\sf r}_-$ to the generator of $\Z/4$ and ${\sf r}_+\mapsto (3,1)\in \Z/4\times \Z/2$. By this description,~$\mathcal{L}$ is 4-periodic: we have isomorphisms $\mathcal{L}^\bullet\cong \mathcal{L}^{\bullet+4k}$ for all~$k$. 

We require the following definition to make contact with ${\rm L}\otimes\C$.

\begin{defn} A section $s\in \Gamma(\Phi_{0,\PA}^{1|2}(X);\mathcal{L}^k)$ is called \emph{supersymmetric} when it is in the image of the restriction map $i^*\colon \Gamma(\Phi_{\epsilon,\PA}^{1|2}(X);\mathcal{L}^k)\to \Gamma(\Phi_{0,\PA}^{1|2}(X);\mathcal{L}^k)$. We denote the set of supersymmetric sections by $\Gamma_{\rm susy}(\Phi_{0,\PA}^{1|2}(X);\mathcal{L}^k)\subset \Gamma(\Phi_{0,\PA}^{1|2}(X);\mathcal{L}^k)$.
\end{defn}

\begin{proof}[Proof of Theorem~\ref{thm1}] We will show there is a natural isomorphism of rings
$$
\Gamma_{\rm susy}(\Phi_{0,\PA}^{1|2}(X);\mathcal{L}^\bullet)\cong \bigoplus_{i\in \Z} \Omega^{\bullet+4i}_{\rm cl}(X)
$$
between supersymmetric sections and 4-periodic closed differential forms. The theorem follows since concordance classes of closed forms are precisely de~Rham cohomology classes, e.g., see the Appendix of~\cite{DBE_WG}. We compute supersymmetric sections as functions on $\R_{>0}\times \pi TX$ that are invariant under the $\R^{1|1}$-action, equivariant for the $T$-action, and in the image of the restriction map that defines supersymmetric sections. 

First we understand the supersymmetric condition, by characterizing the $\R^{1|1}$-action determined by (\ref{eq:PAaction}) on $C^\infty(\R^{1|1}_{>0}\times \pi TX)$.

\begin{lem}
The $\R^{1|1}$-action in (\ref{eq:PAaction}) on $C^\infty(\R^{1|1}_{>0}\times \pi TX)\cong (C^\infty(\R_{>0})[\rho_1])\otimes \Omega^\bullet(X)$ is determined by the formula $\exp(iu Q^2+ \nu_1 Q)$ for the infinitesimal generator
$$
Q:=2i \rho_1 \frac{d}{dr}\otimes {\rm id}-\id\otimes d+i\frac{\rho_1}{r}\otimes {\rm deg}
$$
where $d$ is the de~Rham $d$ and ${\rm deg}$ is the degree endomorphism on differential forms. \end{lem}

This is proved as Lemma~2.11 in~\cite{DBE_WG}, where we use the fact that the formula~(\ref{eq:PAaction}) for the action is identical to the $\R^{1|1}$-action on $\R_{>0}\times \pi TX$ in~\cite{DBE_WG}. 

Returning to the proof of Theorem~\ref{thm1}, it suffices to consider functions $g\in C^\infty(\R_{>0}\times \pi TX)\subset C^\infty(\R^{1|1}_{>0}\times \pi TX)$, i.e., those independent of $\rho$. We write express such a function as
$$
g(r,\rho,x,\psi)=\sum g_l(r)\otimes \alpha_l \quad g_l(r)\in C^\infty(\R_{>0}), \ \alpha_l\in \Omega^\bullet(X).
$$
To be a section of $\mathcal{L}^k$, by definition it must be equivariant for the action of $T$. This in turn demands an equivariance for the $\Z/4$-action generated by ${\sf r}_-$, and invariance under the $\Z/2$-action generated by ${\sf r}_+{\sf r}_-$. This later $\Z/2$-action is trivial, so invariance is automatic. The $\Z/4$ action on~$g$ is by $(i)^{{\rm deg}(\alpha_l)}$, so being a section of $\mathcal{L}^k$ requires that ${\rm deg}(\alpha_l)=k$ mod~4, i.e., we obtain 4-periodic differential forms. 

The remaining condition for $g$ to descend to a section of~$\mathcal{L}$ on $\Gamma_{\rm susy}(\Phi_{0,\PA}^{1|2}(X);\mathcal{L}^k)$ is that it be $Q$-closed. So we compute
$$
Qg=\sum (2i\cdot \rho \frac{dg_l}{dr}\otimes \alpha_l-g_l\otimes d\alpha_l+i\cdot {\rm deg}(\alpha_l)g_l\rho/r\otimes \alpha_l)
$$
where $\alpha_l\in \Omega^k(X)\subset C^\infty(\pi TX)$. We first observe that $Q$-closed functions have $d\alpha_l=0$ for all ${k}$ (e.g., by restricting to the locus $\rho=0$). Furthermore, 
$$
\frac{d g_l}{dr}= \frac{{\rm deg}(\alpha_l)g_l}{2r},
$$
so $g_l(r)=c\cdot r^{{\rm deg}(\alpha_l)/2}$ for some constant $c$, and without loss of generality we may take~$c=1$. Since the dependence on $\R_{>0}$ is therefore completely determined by the degree of the form~$\alpha$, we may identify the function $g$ with a sum of even or odd closed forms, depending on the parity of~$k$. This completes the proof.\ep

The above proof shows that there are various choices of isomorphism between supersymmetric sections and cocycles. In the next section we will use the fixed choice
\beq
\Gamma_{\rm susy}(\Phi_{0,\PA}^{1|2}(X);\mathcal{L}^\bullet)\stackrel{\sim}{\to} \bigoplus_{i\in \Z} \Omega_{\rm cl}^{\bullet+4i}(X),
\quad \quad r^{k/2}\otimes\omega\mapsto \frac{1}{(2\pi)^{k/2}}\omega\label{eq:cochain}
\eeq
for $\omega\in \Omega^k_{\rm cl}(X)$.

\section{Perturbative quantization and the ${\rm L}$-class} 

\subsection{The normal bundles}
Now let $X$ be a Riemannian manifold that we equip with its Levi-Civita connection,~$\nabla$. Given a map $\phi\colon S\times \R^{0|1}\to X$, the pullback bundle $\phi^*TX$ has both a metric and connection. 

\begin{defn} 
Let $\mathcal{N}\Phi_{0,\PA}^{1|2}(X)$ be stack whose $S$-points are triples $(R,\phi,\delta\nu)$ where $(R,\phi)\in \Phi_{0,\PA}^{1|2}(X)(S)$ and $\delta\nu\in \Gamma_0(S\times_{R}\R^{1|2},\phi^*TX)$ where the zero subscript denotes the sections in the orthogonal complement to the constant sections (using the usual pairing on functions on the circle and the Riemannian metric on $X$). Morphisms $(R,\phi,\delta\nu)\to (R',\phi',\delta\nu')$ over $S$ are determined by morphisms $(R,\phi)\to (R',\phi')$ in $\Phi_{0,\PA}^{1|2}(X)(S)$, where $\delta\nu'$ is the pullback of $\delta\nu$ along the map between families of super circles. \end{defn}

\begin{rmk}
The Riemannian exponential map on~$X$ defines an exponential map on sections of~$\mathcal{N}\Phi_{0,\PA}^{1|2}(X)$ with values in~$\Phi^{1|2}(X)$; see~\cite{DBE_WG} Section~4.2 for details.\end{rmk}

We define the linearized classical action as
\beq
\mathcal{S}_{\rm lin}(\phi,\delta\nu)=\int_{S\times_R \R^{1|2}/S} \langle \phi^*\nabla_{D_1} \delta\nu,\phi^*\nabla_{D_2}\delta\nu\rangle,\label{eq:linaction}
\eeq
and we will use the following lemma (proved in \cite{DBE_WG}, Section~4.4) to unravel this definition.

\begin{lem}
Let $\phi_0\colon S\times \R^{0|1}\to X$ be a map with Taylor expansion $\phi_0=x_0+\theta\psi_0$ and $i_0\colon S\to S\times \R^{0|1}$ the inclusion at~$\theta=0$. The $C^\infty(S)[\theta]$-linear map
$$
v\mapsto i_0^*v+ \theta i_0^*((\phi_0^*\nabla_{\partial_\theta}) v),\quad \Gamma(S\times \R^{0|1},\phi_0^*E)\to \Gamma(S,x_0^*E)[\theta]\cong \Gamma(S\times \R^{0|1},p_1^*x_0^*E)
$$
gives an isomorphism of vector bundles over~$S\times \R^{0|1}$. The image of $\phi^*\nabla_{\partial_{\theta}}$ along this inverse of this isomorphism is the operator $\partial_\theta+\theta F(\psi_0,\psi_0)$, where $F(\psi_0,\psi_0)$ denotes the curvature of the connection~$\nabla$ viewed as an ${\rm End}(x_0^*E)$-valued function on~$S$. \label{lem:covarder}
\end{lem}

After pulling back along $\proj$, the above lemma gives 
$$
\phi^*\nabla_{\partial_{\theta_1}}=({\proj}^*\phi_0^*\nabla)_{\partial_{\theta_1}}=\partial_{\theta_1}+\theta_1\mathcal{R}(\psi_1,\psi_1).
$$
The maps $\phi\colon S\times_{R}\R^{1|2}\to X$ are independent of~$t$ and~$\theta_2$, by virtue of factoring through $S\times \R^{0|1}$. Hence, with respect to the isomorphism in Lemma~\ref{lem:covarder} pulled back to $S\times_R\R^{1|2}$, we may identify $\phi^*\nabla_{\partial_t}=\partial_t$ and $\phi^*\nabla_{\partial_{\theta_2}}=\partial_{\theta_2}$, and we identify~$\delta\nu$ with~$a+\theta_1\eta_1+\theta_2\eta_2+\theta_1\theta_2G$ for $a,G$ sections of $TX$ pulled back to $S$ and $\eta_1,\eta_2$ sections of $\pi TX$ pulled back to $S$. Then we compute
\beq
S(\phi,\delta\nu)&=&\int_{S\times_{R}\R^{1|2}/S} \langle \phi^*\nabla_{D_1}\delta\nu,\phi^*\nabla_{D_2} \delta \phi\rangle\nonumber\\
&=&\int_{S\times_{R}\R^{1|2}/S} \langle (\partial_{\theta_1}-i\theta_1\partial_t+\theta_1\mathcal{R}(\psi_1,\psi_1))(a+\theta_1\eta_1+\theta_2\eta_2+\theta_1\theta_2 G),\nonumber \\
 && \phantom{\int_{S\times_{R}\R^{1|2}/S}} (\partial_{\theta_2}-i\theta_2\partial_t)(a+\theta_1\eta_1+\theta_2\eta_2+\theta_1\theta_2 G)\rangle\nonumber\\
 &=&\int_{S^1\times S/S}dt\big( |\dot{a}|^2+i\langle \dot{\eta_2},\eta_2\rangle-i\langle \mathcal{R}(\psi_1,\psi_1)a,\dot{a}\rangle+\langle \mathcal{R}(\psi_1,\psi_1)\eta_2,\eta_2\rangle\nonumber \\
&& \phantom{\int_{S\times_{R}\R^{1|2}/S}}-i\langle \eta_1,\dot{\eta_1}\rangle+\langle G,G\rangle\big). \nonumber\\
&=&\int_{S^1\times S/S} dt\left(-\langle a,D_a a\rangle-i\langle \eta_1,D_{\eta_2}\eta_1\rangle-i\langle \eta_2,D_\eta \eta_2\rangle+\langle G,G\rangle\right)\nonumber\\
&=:&\int_{S^1\times S/S} dt \langle  \delta \nu,\Delta^{1|2}_{X,\PA}\delta \nu\rangle.\nonumber
\eeq
where the sections $a,\eta_1,\eta_2,G$ satisfy the boundary conditions:
$$
a(t)=a(t+r), \quad \eta_1(t)=\eta_1(t+r),\quad \eta_2(t)=-\eta_2(t+r),\quad G(t)=-G(t+r)
$$
and we have separated the even and odd pieces of $\Delta^{1|2}_{X,\PA}$ as 
$$
D_a:={\rm Id}_n\otimes \frac{d^2}{dt^2}-i\mathcal{R}(\psi_1,\psi_1)\otimes \frac{d}{dt},\quad D_{\eta_1}:={\rm Id}_n\otimes \frac{d}{dt}, \quad D_{\eta_2}:={\rm Id}_n\otimes\frac{d}{dt}+i\mathcal{R}(\psi_1,\psi_1)\otimes {\rm id}, 
$$
where ${\rm Id}_n$ denotes the identity operator on sections of the pullback of $TX$ and $\pi TX$. Following the standard rules of functional integration in this easy case, we define the $\zeta$-super determinant as
$$
{\rm sdet}_\zeta(\Delta^{1|2}_{X,\PA}):=\frac{{\rm Pf}_\zeta(D_{\eta_1}){\rm Pf}_\zeta(D_{\eta_2})}{{\rm det}(D_a)^{1/2}}.
$$

\subsection{The ${\rm L}$-class}

The ${\rm L}$-class is associated to the characteristic series 
$$
\frac{x/2}{\tanh(x/2)}=\frac{x\cosh(x/2)}{2\sinh(x/2)}
$$
From the product formulas
$$
\sinh(x/2)=(x/2)\prod_{k=1}^\infty \left(1+\left(\frac{x}{2\pi k}\right)^2\right),\quad \cosh(x/2)=\prod_{k=1}^\infty\left(1+\left(\frac{x}{2\pi(k-1/2)}\right)^2\right)
$$
we derive 
$$
\frac{\sinh(x/2)}{x/2}=\exp\Big(-\sum_{k=1}^\infty \frac{x^{2k}}{2k(2\pi i)^{2k}} 2\zeta(2k)\Big)  ,\ \cosh(x)=\exp\Big(-\sum_{k=1}^\infty \Big(\frac{x^{2k}}{2k(2\pi i)^{2k} } 2\sum_{n=1}^\infty \frac{1}{(n-1/2)^{2k}}\Big)\Big),
$$
where $\zeta$ denotes the Riemann $\zeta$-function. This will give us a convenient form for the ${\rm L}$-class. 

\begin{proof}[Proof of Theorem~\ref{thm2}]
To compute the $\zeta$-super determinant we require two facts. First, the $\zeta$-regularized product associated to the sequence $\{\frac{(2\pi k)^n}{r^n}\}_{k=1}^\infty$ is~$r^{n/2}$; e.g., see Example~2 of~\cite{zetaprod}. Second, if $A$ is determinant class, then ${\rm det}_\zeta(AD)={\rm det}_{\rm Fr}(A)\cdot {\rm det}_\zeta(D)$ where ${\rm det}_{\rm Fr}$ denotes the Fredholm determinant, and similarly for $\zeta$-regularized Pfaffians. We claim
\beq
{\rm det}_\zeta(D_a)&=&{\rm det}_\zeta\left(\frac{d^2}{dt^2}\otimes {\rm Id}_{TX}\right){\rm det}_{\rm Fr}\left({\rm Id}-i\mathcal{R}(\psi_1,\psi_1)\otimes \left(\frac{d}{dt}\right)^{-1}\right)\nonumber\\
{\rm pf}_\zeta(D_{\eta_2})&=&{\rm pf}_\zeta\left(\frac{d}{dt}\otimes {\rm Id}_{TX}\right){\rm pf}_{\rm Fr}\left({\rm Id}+i\mathcal{R}(\psi_1,\psi_1)\otimes \left(\frac{d}{dt}\right)^{-1}\right).\nonumber
\eeq
To verify that the claimed operators are determinant class, first we use functional calculus to write
\beq
{\rm det}_{\rm Fr}\left({\rm Id}-i\mathcal{R}(\psi_1,\psi_1)\otimes \left(\frac{d}{dt}\right)^{-1}\right)&=&\exp\left(-\sum_{k=0}^\infty \frac{{\rm Tr}((i\mathcal{R}(\psi_1,\psi_1))^{k}\otimes (d/dt)^{-k})}{k} \right)\nonumber\\
{\rm pf}_{\rm Fr}\left({\rm Id}+i\mathcal{R}(\psi_1,\psi_1)\otimes \left(\frac{d}{dt}\right)^{-1}\right)&=&\exp\left(\frac{1}{2}\sum_{k=0}^\infty (-1)^{k+1}\frac{{\rm Tr}((i\mathcal{R}(\psi_1,\psi_1))^k\otimes (d/dt)^{-k})}{k}\right).\nonumber
\eeq
The traces of odd powers of the curvature vanish and the operators~$(d/dt)^{-2k}$ are trace class for~$k>0$, so the claimed operators are indeed determinant class. 

It remains to compute the determinant. For periodic (respectively, antiperiodic) boundary conditions, a basis of sections is given by $\{e^{2\pi i l/r}\}_{l\in \Z_*}$ where $\Z_*=\Z-\{0\}$ (respectively, $\{e^{2\pi i (l-1/2)/r}\}_{l\in \Z}$). 
We have in these respective cases
\beq
\begin{array}{ccccccc}
{\rm Tr}((d/dt)^{-2k})&=&2\sum_{l=1}^\infty \frac{r^{2k}}{(2\pi i l)^{2k}}&=&\frac{2r^{2k} }{(2\pi i)^{2k}} \zeta(2k), & {\rm (periodic)}\nonumber\\
 {\rm Tr}((d/dt)^{-2k})&=&2\sum_{l=1}^\infty \frac{r^{2k}}{(2\pi i (l-1/2))^{2k}}&=&\frac{2r^{2k} }{(2\pi i)^{2k}} \sum_{l=1}^\infty \frac{1}{(l-1/2)^{2k}} & {\rm (antiperiodic).} \end{array}\nonumber
\eeq

Turning attention to the curvature contribution, for $F$ the curvature 2-form of a real vector bundle $V$ with connection, the differential-form valued Pontryagin character of $V$ is defined by
\beq
{\rm Tr}(F^{2k})=(2k)!(2\pi i)^{2k} {\rm ch}_{2k}(V\otimes \C)=(2k)!(2\pi i)^{2k} {\rm ph}_{k}(V),\label{eq:formChern}
\eeq
where ${\rm ch}_{2k}$ denotes the $4k^{\rm th}$ component of the Chern character, and ${\rm ph}_{k}$ denotes the $4k^{\rm th}$ component of the Pontryagin character. In our cochain model (\ref{eq:cochain}) we have
$$
(ir)^{2k}(1/2){\rm Tr}\left(\mathcal{R}\right)^{2k}=(2k)! {\rm ph}_{k}(TX),
$$
where (in an abuse of notation) the right hand side above denotes the $k^{\rm th}$ component of the Pontryagin character as a function on~$\Phi^{1|2}_{0,\PA}(X)$. 

Putting this together, we obtain 
\beq
{\rm det}_\zeta(D_a)&=& r^{2n}\exp\left(-\sum_{k=1}^\infty \frac{(2k)!{\rm ph}_k(TX)}{(2k)(2\pi i)^{2k}} 2\zeta(2k)\right)\nonumber\\
{\rm pf}_\zeta(D_{\eta_2})&=&r^{n/2}\exp\left(-\sum_{k=1}^\infty \frac{(2k)!{\rm ph}_k(TX)}{2k(2\pi i)^{2k}} 2 \sum_{l=1}^\infty \frac{1}{(l-1/2)^{2k}}\right)\nonumber\\
{\rm pf}_\zeta(D_{\eta_1})&=&r^{n/2}\nonumber
\eeq
and hence the asserted ratio is a representative for the ${\rm L}$-class as a function on~$\Phi^{1|2}_{0,\PA}(X)$. It is supersymmetric because the ${\rm ph}_k(TX)$ are. 

Lastly, there is a canonical relative volume form on the fibers $\int_X\colon \Phi_{0,\PA}^{1|2}(X)\to \Phi_{0,\PA}^{1|2}(\pt)$ gotten from multiplication by $r^{n/2}$ and integration of differential forms. If we modify this volume form by the ${\rm L}$-class as constructed by the above super determinant, we obtain the map
\beq
 \Gamma_{\rm susy}(\Phi_{0,\PA}^{1|2}(X);\mathcal{L}^k) &\stackrel{\pi_!}{\to} & \Gamma_{\rm susy}(\Phi_{0,\PA}^{1|2}(\pt);\mathcal{L}^{k-n})\nonumber\\
 s&\mapsto& \int_Xs\cdot {\rm sdet}_\zeta(\Delta^{1|2}_{X,\PA}), \nonumber
\eeq
which we identify with the ${\rm MSO}$-orientation of ${\rm L}$-theory tensored with~$\C$.
\ep

\bibliographystyle{amsalpha}
\bibliography{references}

\end{document}